\documentclass{elsart2}

\usepackage{amssymb}

\renewcommand{\(}{\left(}
\renewcommand{\)}{\right)}

\newcommand{\im}{{\rm im}}

\newcommand{\K}{\wh{K}}
\newcommand{\got}{\mathfrak}

\newcommand{\X}{X_1,\ldots ,X_n}

\newcommand{\Ym}{Y_1,\ldots ,Y_m}
\newcommand{\Y}{Y_1,\ldots ,Y_n}

\renewcommand{\v}{\wh{v}}

\newcommand{\wh}{\widehat}
\newcommand{\om}{\omega}
\newcommand{\Om}{\Omega}
\newcommand{\al}{\alpha}

\newcommand{\si}{\sigma}

\newcommand{\De}{\Delta}

\newcommand{\be}{\beta}
\newcommand{\vp}{\varphi}

\newcommand{\R}{R_{\v}}
\renewcommand{\r}{R_{v}}
\newcommand{\M}{{\got m}_{\v}}
\newcommand{\m}{{\got m}_v}

\newcommand{\Kn}{{\widehat K_n}}
\newcommand{\D}{\Delta_{\v}}
\renewcommand{\d}{\Delta_v}

\newcommand{\Phii}{\Phi^{-1}}
\newcommand{\pst}{\Phi_{\sigma ,\theta}}
\newcommand{\psti}{\(\pst\)^{-1}}
\newcommand{\lcor}{{\rm [\kern - 1.8pt [}}
\newcommand{\rcor}{{\rm ]\kern - 1.8pt ]}}
\renewcommand{\[}{\left[}
\renewcommand{\]}{\right]}
\newcommand{\llcor}{\[\kern -2.5pt\[}
\newcommand{\rrcor}{\]\kern -2.5pt\]}
\newcommand{\Dt}{\D\lcor t\rcor}

\newcommand{\fb}{\overline{f}}

\newcommand{\ZZ}{{\mathbb{Z}}}

\newcommand{\CC}{{\mathbb{C}}}
\newcommand{\FF}{{\mathbb{F}}}

\newcommand{\LL}{{\mathbb{L}}}
\newcommand{\lto}{\longrightarrow}
\newcommand{\lmapsto}{\longmapsto}

\begin{document}

\sloppy

\begin{frontmatter}

\title{Rank one discrete valuations of $k((X_1,\ldots ,X_n))$}
\author{Miguel \'{A}ngel Olalla--Acosta\thanksref{MAO}}
\address{Facultad de Matem\'aticas. Apdo. 1160. E-41080 SEVILLA (SPAIN)}
\thanks[MAO]{Partially supported by Junta de Andaluc\'{\i}a, Ayuda a
grupos FQM 218}
\ead{maolalla@us.es}

\begin{abstract}
In this paper we study the rank one discrete valuations of
$k((X_1,\ldots ,X_n))$ whose center in $k\lcor\X\rcor$ is the maximal
ideal $(\X )$. In sections 2 to 6 we give a construction of a
system of parametric equations describing such valuations. This amounts to
finding a parameter and a field of coefficients. We devote section 2
to finding an element of value 1, that is, a parameter. The field of
coefficients is the residue field of the valuation, and it is
given in section 5.

The constructions given in these sections are not effective in the
general case, because we need either to use the Zorn's lemma or to know
explicitly a section $\sigma$ of the natural homomorphism $R_v\to\d$
between the ring and the residue field of the valuation $v$.

However, as a consequence of this construction, in section 7, we
prove that $k((\X ))$ can be embedded into a field $L((\Y ))$, where
the {\em ``extended valuation'' is as close as possible to the usual
  order function}. 
\end{abstract}

\begin{keyword}
Valuation theory, local uniformization, formal power rings, completions.
\MSC 13F25 \sep 13F30 \sep 14B05 \sep 16W60
\end{keyword}

\end{frontmatter}

\section{Terminology and preliminaries}

Let $k$ be an algebraically closed field of characteristic 0,
$R_n=k\lcor\X\rcor$, $M_n=(\X )$ its maximal ideal and $K_n=k ((\X
))$ its quotient field. Let $v$ be a rank-one discrete valuation
of $K_n\vert k$, $\r$ the valuation ring, $\m$ the maximal ideal
and $\d$ the residue field of $v$. The center of $v$ in $R_n$ is
$\m\cap R_n$. Throughout this paper ``discrete valuation of
$K_n\vert k$" means ``rank-one discrete valuation of $K_n\vert k$
whose center in $R_n$ is the maximal ideal $M_n$". The dimension of
$v$ is the transcendence degree of $\d$ over $k$. In order to simplify
the redaction we shall assume, without loss of generality, that the
group of $v$ is $\ZZ$.

Let $\Kn$ be the completion of $K_n$ with respect to $v$, $\v$ the
extension of $v$ to $\Kn$, $\R$, $\M$ and $\D$ the ring, maximal
ideal and the residue field of $\v$, respectively (see
\cite{Ser1} for more details). We know that $\d$ and $\D$ are isomorphic
(\cite{Krull}). Let $\si :\D\to\R$ be a $k-$section of the natural
homomorphism $\R\to\D$, $\theta\in\R$ an element of value 1 and
$t$ an indeterminate. We consider the $k-$isomorphism $$\Phi =
\pst :\Dt\to\R$$ given by $$\Phi\(\sum\al_it^i\) =\sum\si
(\al_i)\theta^i,$$ and denote also by $\Phi$ its extension to the
quotient fields. We have a $k-$isomorphism $\Phii$ which, when
composed with the usual order function on $\D ((t))$, gives the
valuation $\v$. This is the situation we will consider throughout
this paper, and we will freely use it without new explicit
references.

We shall use two basic transformations in order to find
an element of value 1 and construct the residue field:

\begin{enumerate}

\item {\em Monoidal transformation:} $$ \begin{array}{rcl} k\lcor\X\rcor
& \lto & k\lcor\Y\rcor \\ X_1 & \lmapsto & Y_1 \\ X_2 & \lmapsto &
Y_1Y_2 \\ X_i & \lmapsto & Y_i,\ i=3,\ldots ,n.
\end{array}
$$ with $v(X_2)>v(X_1)$.

\item {\em Change of coordinates:}$$
\begin{array}{rcl}
k\lcor\X\rcor & \lto & L\lcor\Y\rcor \\ X_1 & \lmapsto & Y_1 \\ X_i
& \lmapsto & Y_i+c_iY_1,\ i=2,\ldots ,n,
\end{array}
$$ where $c_i\in\R\setminus\M$ and $L$ is an extension field of $k$.

\end{enumerate}

For both transformations we have the following facts:

\begin{enumerate}

\item[(a)] The transformations are one to one: In the case of the monoidal
transformations this property is well known. In the other case it is a
consequence of \cite{ZSII} (corollary 2, page 137).


\item[(b)] New variables $Y_i$ lie in $\R$, so we can put $\psti
(Y_i)=\sum a_{i,j}t^j$.

\item[(c)] Let $\varphi :K_n\to\d ((t))$ be the restriction of $\psti$ to
$K_n$. Let us denote by $\varphi ':L_n=L((\Y ))\to\d ((t))$ the natural
extension of $\varphi$ to the field $L_n$. Then $v=\nu_t\circ\varphi
'_{\vert K_n}$, with $\nu_t$ the usual order function over $\d
((t))$. Therefore, if $\varphi '$ is injective we can extend the
valuation $v$ to the field $L_n$ and the extension is
$v'=\nu_t\circ\varphi '$.

\end{enumerate}

From now on transformation will mean monoidal
transformation, change of coordinates, variables interchanges or
finite compositions of these.


\section{Construction of an element of value 1}

Remember that we are assuming that the group of $v$ is $\ZZ$, so there
exists an element $u\in K$ such that $v(u)=1$.

\begin{lem}{\label{lema1}}
Let $\al_i=v(X_i)$ for all $i=1,\ldots, n$. By a finite number of monoidal
transformations we can find $n$ elements $Y_1,\ldots ,Y_n\in\Kn$ such
that $v (Y_i)=\al =\gcd\{\al_1,\ldots ,\al_n\}$.
\end{lem}

\begin{pf}
We can suppose that $v(X_1)=\al_1=\min\{\al_i\vert 1\leq i\leq
n\}$ and consider the following two steps:

{\bf Step 1.-} If there exists $n_i\in\ZZ$ such
that $\al_i=n_i\al_1$ for all $i= 2,\ldots , n$, then for each $i$ we
apply $n_i-1$ monoidal
transformations
$$ \begin{array}{rcl} k\lcor\X\rcor & \lto &
k\lcor\Y\rcor \\
X_i & \lmapsto & Y_1Y_i\\
X_j & \lmapsto & Y_j,\ j\ne i.
\end{array}
$$ Trivially $v (Y_i)=\al_1$ for all $i=1,\ldots ,n$.

{\bf Step 2.-} Assume there exists $i$, with $2\leq i\leq n$, such that
$v(X_1)=\al_1$ does not divide to $v(X_i)=\al_i$. We can suppose
that $i=2$ with no loss of generality and then $\al_2 = q\al_1+r$. So
we apply $q$ times the monoidal transformation
$$
\begin{array}{rcl} k\lcor\X\rcor & \lto & k\lcor\Y\rcor \\
X_2 & \lmapsto & Y_1Y_2 \\
X_i & \lmapsto & Y_i,\
i\ne 2
\end{array}
$$ to obtain a new ring $k \lcor\Y\rcor$ where $v (Y_2) =r>0$ and
$Y_2$ is the element of minimum value.

As the values of the variables are greater than zero, in a finite
number of steps 2 we come to the situation of step 1. In fact,
this algorithm is equivalent to the ``euclidean algorithm" to
compute the greatest common divisor of $\al_1,\ldots ,\al_n$. \qed
\end{pf}

\begin{thm}
We can construct an element of value 1 applying a finite number
of monoidal transformations and changes of coordinates.
\end{thm}

\begin{pf}
We call $Y_{1,r},\ldots ,Y_{n,r}$ the elements found after $r$
transformations.

We can suppose that we have applied the previous lemma
to obtain $Y_{1,1},\ldots ,Y_{n,1}$ such that $v(Y_{i,1})=\al$ for all
$i=1,\ldots ,n$. Let us prove that there exists
$c_i\in\R\setminus\M$ for each $i=2,\ldots ,n$ such that $\v
(Y_{i,1}-c_iY_{1,1}) >\al$. We can take
$$\psti (Y_{i,1})=\sum_{j\geq\al} a_{i,j}t^j=\om_i(t),\
a_{i,j} \in\D,\ a_{i,\al}\ne 0,$$
and so it suffices taking $b_i=a_{i,\al}/a_{1,\al }$ and $c_i=\si (b_i)$.

The following two steps defines a procedure to obtain an element of value 1:

{\bf Step 1.-} We apply the coordinate change
$$
\begin{array}{rcl}
k\lcor Y_{1,1},\ldots ,Y_{n,1}\rcor & \lto & L\lcor Y_{1,2},\ldots
,Y_{n,2}\rcor \\
Y_{1,1} & \lmapsto & Y_{1,2} \\
Y_{i,1} & \lmapsto & Y_{i,2}+ c_iY_{1,2},\ i=2,\ldots ,n.
\end{array}
$$
With this transformation the values of the new variables are not
equal to $\v (Y_{1,2})$.

{\bf Step 2.-} We apply lemma \ref{lema1} to equalize the values of
elements and go to step 1. Obviously, the
minimum of the values of the elements does not increase, because the
greater common divisor of the values does not exceed the minimum of
the values. Moreover the first variable does not change.

If we obtain an element of value 1 then we are finished.

We have to show that the procedure produces an element of value
1 in a finite number of transformations. The only way for the process
to be infinite is that, in step 2, the minimum of the values of
the elements does not decrease. This means that, in step 1, the value
of the first variable divides the values of the new variables.

The composition of steps 1 and 2 is the transformation
$$
\begin{array}{rcl}
k\lcor Y_{1,r},\ldots ,Y_{n,r}\rcor & \lto & L\lcor Y_{1,r+1}\ldots
,Y_{n,r+1}\rcor \\
Y_{1,r} & \lmapsto & Y_{1,r+1} \\
Y_{i,r} & \lmapsto & Y_{i,r+1}+c_iY_{1,r+1}^{m_i},\ i=2,\ldots ,n.
\end{array}
$$

If we use steps 1 and 2 infinitely many times, we have an infinite
sequence of transformations
$$
\begin{array}{rcl}
k\lcor\Y\rcor & \lto & L\lcor Y_{1,j},\ldots ,Y_{n,j}\rcor \\
Y_1 & \lmapsto & Y_{1,j}
\\ Y_i & \lmapsto & Y_{i,j}+\sum_{k=1}^j c_{i,k}Y_{1,j}^{m_{i,k}},\
i=2,\ldots ,n.
\end{array}
$$
Then we can obtain an infinite sequence of variables
$$
\begin{array}{rcl}
Y_{1,j} & = & Y_{1,j} \\
Y_{i,j} & = & Y_i-\sum_{k=1}^j c_{i,k}Y_1^{m_{i,k}},\
i=2,\ldots ,n,
\end{array}
$$
with $\v (Y_{i,j})>\v (Y_{i,j-1})$ for all $i,\ j$. So any
sequence of partial sums of the series
$$Y_i-\sum_{k=1}^{\infty} c_{i,k}Y_1^{m_{i,k}},\ \forall
i=2,\ldots ,n$$
have strictly increasing values. Then these series
converge to zero in $\R$, so
$$Y_i=\sum_{k=1}^{\infty}
c_{i,k}Y_1^{m_{i,k}},\ \forall i=2,\ldots ,n.$$
Let $f(Y_1,\ldots ,Y_n)\in K_n$, then
$$v(f)=\v\( f\(
Y_1,\sum_{k=1}^{\infty} c_{2,k}Y_1^{m_{2,k}}, \ldots
,\sum_{k=1}^{\infty} c_{n,k}Y_1^{m_{n,k}}\)\) =m\cdot
v(Y_1).$$
In this situation, the  group of $v$ is
$v(Y_1)\cdot\ZZ$ (see \cite{Bri2}) but as the group is assumed to be
$\ZZ$, $\v (Y_1)=1$. \qed
\end{pf}

\begin{exmp}
Let us consider the embedding
$$
\begin{array}{rcl}
\Psi :  \CC\lcor X_1,X_2,X_3\rcor & \lto &
        \CC (T_2,T_3)\lcor t\rcor \\
X_1 & \lmapsto & t^2 \\
X_2 & \lmapsto & T_2t^4+T_2t^6 \\
X_3 & \lmapsto & T_2t^2+T_3t^5
\end{array}
$$
with $t$, $T_2$ and $T_3$ variables over $\CC$. We are going to
denote its extension to the quotient fields by $\Psi$ as well. The
composition of this injective homomorphism with the order function
in $t$ gives a discrete valuation of $\CC
((X_1,X_2,X_3))\vert\CC$, $v=\nu_t\circ\Psi$. If we apply the
procedure given in this section we construct the following
element of value 1: $$\frac{X_3- c_2X_1}{X_1^2},$$ where
$c_2\in\R\setminus\M$ such that $\Psi (c_2)=T_2+t\cdot f,$ with $f\in\CC
(T_2,T_3)\lcor t\rcor$. In this case we can take 
$$c_2=\frac{X_2}{X_1^2+X_1^3}.$$
\end{exmp}

\begin{rem}
We need to know some elements $c_i\in\R\setminus\M$ (or $b_i\in\D$
and $\si :\D\to\R$) such that $\v (Y_{i,1}-c_iY_{1,1})>\al$ for each
$i=2,\ldots ,n$ in order to apply the procedure described in the proof of
theorem 2. Let $\Delta$ be a field, if the valuation is given as a
composition of an injective homomorphism
$$
\begin{array}{rcl}
\Psi :k\lcor\X\rcor& \lto & \Delta\lcor t\rcor \\
X_i & \lmapsto & \sum_{j\geq 1}a_{i,j}t^j
\end{array}
$$
with the usual order function of $\Delta ((t))$, $v=\nu_t\circ\Psi$, then we
can find the $c_i$'s using the coefficients $a_{i,j}\in\Delta$ of $\Psi
(X_i)$.
\end{rem}

\section{Transcendental and algebraic elements of $\d$}

In the following sections we give a procedure to construct the
residue field $\d$ of a discrete valuation of $K_n\vert k$, as a
transcendental extension of $k$.

Before the description of the procedure we have to do the following remark
about the $k-$section $\si$.

\begin{rem}{\label{remark4}}
We are going to check all the variables searching those residues which 
generate the extension $k\subset\D$.
Hence we will have to move between $\R$ and $\D$ by the $k-$section
$\si$ and the natural homomorphism $\D\to\R$. We can do the following
considerations:

\noindent {\bf 1)} Let us consider the diagram
\begin{center}
\unitlength=0.75mm  
\linethickness{0.4pt}
\begin{picture}(49.00,60.00)
\vspace{.5cm}
\put(10.00,4.00){\makebox(0,0)[cc]{$k$}}
\put(29.50,4.33){\vector(1,0){14.00}}
\put(29.50,4.33){\vector(-1,0){14.00}}
\put(28.00,6.33){\makebox(0,0)[cc]{{\footnotesize $id$}}}
\put(49.00,4.00){\makebox(0,0)[cc]{$k$}}
\put(10.00,7.00){\vector(0,1){15.00}}
\put(49.00,7.00){\vector(0,1){15.00}}
\put(10.00,27.67){\makebox(0,0)[cc]{$\FF$}}
\put(10.00,32.33){\vector(0,1){15.00}}
\put(49.00,32.33){\vector(0,1){15.00}}
\put(10.00,52.67){\makebox(0,0)[cc]{$\R$}}
\put(49.00,52.67){\makebox(0,0)[cc]{$\D$}}
\put(16.00,26.00){\vector(1,0){28.00}}
\put(16.00,51.00){\vector(1,0){28.00}}
\put(44.00,28.00){\vector(-1,0){28.00}}
\put(44.00,53.00){\vector(-1,0){28.00}}
\put(28.00,24.00){\makebox(0,0)[cc]{{\footnotesize $\vp$}}}
\put(28.00,49.00){\makebox(0,0)[cc]{{\footnotesize $\vp$}}}
\put(28.00,30.00){\makebox(0,0)[cc]{{\footnotesize $\si$}}}
\put(28.00,55.00){\makebox(0,0)[cc]{{\footnotesize $\si$}}}
\put(49.00,27.67){\makebox(0,0)[cc]{$\FF'$}}
\end{picture}
\end{center}
where $\FF$ and $\FF '$ are subfields of $\R$ and $\D$
respectively. Let $\om\in\R$ an element such that $\v (\om )=0$.
The question is: if $\om +\M$ is transcendental over $\FF '$, is
$\si (\om +\M )$ transcendental over $\FF$? What happens in the
algebraic case?

So we suppose $\om +\M$ to be transcendental over $\FF '$. Let
$f(X)\in\FF [X]$ be a non-zero polynomial. Let us put
$$f(X)=\sum_{i=0}^n\si (a_i')X^i,\ a_i'\in\FF '.$$
Then
$$f(\si
(\om +\M ))=\sum_{i=0}^n\si (a_i')\si (\om +\M
)^i=\si\(\sum_{i=0}^na_i'(\om +\M )^i\)\ne 0$$
because $\om +\M$
is transcendental over $\FF '$. So we have proved that $\si (\om + \M)$
is transcendental over $\FF$ if $\om +\M$ is transcendental over
$\FF '$

\noindent {\bf 2)} In the algebraic case let us consider the next
diagram: \begin{center}
\unitlength=0.75mm
\linethickness{0.4pt}
\begin{picture}(49.00,59.67)
\vspace{.5cm} \put(10.00,4.00){\makebox(0,0)[cc]{$k$}}
\put(29.50,4.33){\vector(1,0){14.00}}
\put(29.50,4.33){\vector(-1,0){14.00}}
\put(28.00,6.33){\makebox(0,0)[cc]{{\footnotesize $id$}}}
\put(49.00,4.00){\makebox(0,0)[cc]{$k$}}
\put(10.00,7.00){\vector(0,1){15.00}}
\put(49.00,7.00){\vector(0,1){15.00}}
\put(10.00,27.67){\makebox(0,0)[cc]{$\FF$}}
\put(10.00,32.33){\vector(0,1){15.00}}
\put(49.00,32.33){\vector(0,1){15.00}}
\put(10.00,52.67){\makebox(0,0)[cc]{$\R$}}
\put(49.00,52.67){\makebox(0,0)[cc]{$\D$}}
\put(16.00,26.00){\vector(1,0){28.00}}
\put(44.00,28.00){\vector(-1,0){28.00}}
\put(28.00,24.00){\makebox(0,0)[cc]{{\footnotesize $\vp$}}}
\put(28.00,30.00){\makebox(0,0)[cc]{{\footnotesize $\si$}}}
\put(49.00,27.67){\makebox(0,0)[cc]{$\FF'$}}
\end{picture}
\end{center}
Let $\al +\M\in\D$ be an algebraic element over $\FF '$, with $\v
(\al )=0$ (i.e. $\al +\M\ne 0$). Let
$$\fb (X) =
X^m+\be_1X^{m-1}+\cdots +\be_m\in\FF '[X]$$
be its minimal polynomial over $\FF '$. Let us take the polynomial
$$f(X) =
X^m+b_1X^{m-1}+\cdots +b_m\in\FF[X],\mbox{ with } b_i=\si
(\be_i).$$
By Hensel's Lemma (\cite{ZSII}, corollary 1, page 279) we know
that there exists $a\in\R$ such that $a$ is a simple root of $f(X)$ y
$\varphi (a)=\al +\M$. As $\vp\si =id$, $f(X)$ is the minimal
polynomial of $a$, so we can extend $\si :\FF '[\al +\M ]\to\FF
[a]$. Then we have
\vspace{0.5cm}
\begin{center}
\unitlength=0.75mm
\linethickness{0.4pt}
\begin{picture}(49.00,77.34)
\vspace{.5cm} \put(10.00,4.00){\makebox(0,0)[cc]{$k$}}
\put(29.50,4.33){\vector(1,0){14.00}}
\put(29.50,4.33){\vector(-1,0){14.00}}
\put(28.00,6.33){\makebox(0,0)[cc]{{\footnotesize $id$}}}
\put(49.00,4.00){\makebox(0,0)[cc]{$k$}}
\put(10.00,7.00){\vector(0,1){15.00}}
\put(49.00,7.00){\vector(0,1){15.00}}
\put(10.00,27.67){\makebox(0,0)[cc]{$\FF$}}
\put(10.00,32.33){\vector(0,1){15.00}}
\put(49.00,32.33){\vector(0,1){15.00}}
\put(10.00,51.67){\makebox(0,0)[cc]{$\FF(a)$}}
\put(49.00,51.67){\makebox(0,0)[cc]{$\FF'(\al +\M )$}}
\put(10.00,57.00){\vector(0,1){15.00}}
\put(49.00,57.00){\vector(0,1){15.00}}
\put(10.00,77.34){\makebox(0,0)[cc]{$\R$}}
\put(49.00,77.34){\makebox(0,0)[cc]{$\D$}}
\put(16.00,26.00){\vector(1,0){28.00}}
\put(44.00,28.00){\vector(-1,0){28.00}}
\put(28.00,24.00){\makebox(0,0)[cc]{{\footnotesize $\vp$}}}
\put(28.00,30.00){\makebox(0,0)[cc]{{\footnotesize $\si$}}}
\put(49.00,27.67){\makebox(0,0)[cc]{$\FF'$}}
\put(20.00,50.67){\vector(1,0){15.00}}
\put(35.00,52.67){\vector(-1,0){15.00}}
\put(28.00,48.67){\makebox(0,0)[cc]{{\footnotesize $\vp$}}}
\put(28.00,54.67){\makebox(0,0)[cc]{{\footnotesize $\si$}}}
\end{picture}
\end{center}

Let us consider the set
$$\Om =\{
(\FF_1,\si_1)\vert\FF_1\supset\FF\mbox{ and }\si_1\mbox{ extends
}\si\}$$
partially ordered by
$$(\FF_1,\si_1)<(\FF_2,\si_2)\iff\FF_1\subset\FF_2\mbox{ and  }
\si_{2\vert\FF_1}=\si_1.$$
By Zorn's Lemma there exists a maximal
element $(\LL ,\si ')\in\Om$, and again by Hensel's
Lemma (\cite{ZSII}, corollary 2, page 280) we have $\vp (\LL
)=\D$. So we can extend $\si$ to a $k-$section $\si '$ of $\vp$ in
such a way that $a=\si '(\al +\M)$ is an algebraic element over
$\FF$.

\noindent {\bf 3)} Hence we have showed that if $\om +\M\in\D$, $\v
(\om )=0$, is a transcendental (resp. algebraic) element over $\FF '$,
there exists a $k-$section of $\vp$ which extends $\si$ and $\si (\om
+\M )$ is transcendental (resp. algebraic) over $\FF$.
\begin{center}
\unitlength=0.75mm
\linethickness{0.4pt}
\begin{picture}(49.00,54.00)
\vspace{.5cm} \put(10.00,4.00){\makebox(0,0)[cc]{$k$}}
\put(29.50,4.33){\vector(1,0){14.00}}
\put(29.50,4.33){\vector(-1,0){14.00}}
\put(28.00,6.33){\makebox(0,0)[cc]{{\footnotesize $id$}}}
\put(49.00,4.00){\makebox(0,0)[cc]{$k$}}
\put(10.00,7.00){\vector(0,1){15.00}}
\put(49.00,7.00){\vector(0,1){15.00}}
\put(10.00,27.67){\makebox(0,0)[cc]{$\FF$}}
\put(10.00,31.33){\vector(0,1){15.00}}
\put(49.00,31.33){\vector(0,1){15.00}}
\put(10.00,51.67){\makebox(0,0)[cc]{$\R$}}
\put(49.00,51.67){\makebox(0,0)[cc]{$\D$}}
\put(16.00,26.00){\vector(1,0){28.00}}
\put(44.00,28.00){\vector(-1,0){28.00}}
\put(28.00,24.00){\makebox(0,0)[cc]{{\footnotesize $\vp$}}}
\put(28.00,30.50){\makebox(0,0)[cc]{{\footnotesize $\si$}}}
\put(49.00,27.67){\makebox(0,0)[cc]{$\FF'$}}
\put(16.00,50.00){\vector(1,0){28.00}}
\put(44.00,52.00){\vector(-1,0){28.00}}
\put(28.00,48.00){\makebox(0,0)[cc]{{\footnotesize $\vp$}}}
\put(28.00,54.50){\makebox(0,0)[cc]{{\footnotesize $\si '$}}}
\end{picture}
\end{center}
\end{rem}

\section{A first transcendental residue.}

We devote this section to finding a first transcendental residue of $\d$
over $k$. Note that this preliminary transformations construct the
residue field in the case $n=2$.

\begin{lem}{\label{lema6}}
There exists a finite number of monoidal transformations and
changes of coordinates that constructs $n$ elements $Y_1,\ldots ,Y_n$ such that
$v(Y_i)=v(Y_1)=\al$ and the residue $Y_2/Y_1+\M$ is not in $k$.
\end{lem}

\begin{pf}
We can suppose that we have applied lemma
\ref{lema1} to obtain $Y_1.\ldots ,Y_n$ such that $v(Y_i)=\al$ for
all $i=1, \ldots ,n$.

In this situation $v(Y_i/Y_j)=0$, so $0\ne
(Y_i/Y_j)+\m\in\d$. If this residue lies in $k$ then there exists
$a_{i,1}\in k$ such that 
$$\frac{Y_i}{Y_j}+\m = a_{i,1}+\m ,$$
so
$$\frac{Y_i}{Y_j}-a_{i,1}=\frac{Y_i-a_{i,1}Y_j}{Y_j}\in\m ,$$
and then
$$v\(\frac{Y_i-a_{i,1}Y_j}{Y_j}\) >0.$$
So we have $v(Y_i-a_{i,1}Y_j)=\al_1>\al$. If $\al$ divides to $\al_1$
then $\al_1=r_1\al$ with $r_1\geq 2$ and
$$v\(\frac{Y_i-a_{i,1}Y_j}{Y_j^{r_1}}\) =0.$$
If the residue of this element lies too in $k$, then exist $a_{i,r_1}\in
k$ such that
$$v(Y_i-a_{i,1}Y_j-a_{i,r_1}Y_j^{r_1})=\al_2>\al_1.$$
If $\al$ divides to $\al_2$ then $\al_2=r_2\al$ with $r_2>r_1$ and we
can repeat this operation.

The above procedure is finite for some pair
$(i,j)$. We know (\cite{Bri2}) that any discrete valuation of
$k((X_1,X_2))$ has dimension 1, so the restriction, $v'$, of our valuation $v$
to the field $k((X_1,X_2))$ is a valuation with dimension 1, and the
dimension of $v$ is greater or equal than 1, because a
transcendental residue of $v'$ over $k$ is a transcendental residue of
$v$ too. If the procedure never ends for all $(i,j)$ then all the
residues of $v$ are in $k$, so the dimension of $v$ is 0 and there is
a contradiction. So we can suppose that the above procedure ends
for $(1,2)$ by reordering the variables if necessary.

Hence there exists a first transcendental residue. We can apply the above
procedure to the variables $Y_1,Y_2$, and so we have the transformations:
$$Z_i=Y_i,\ i\ne 2$$
$$Z_2= Y_2-\sum_{i=1}^{s_2}a_{2,i}Y_1^i,$$
such that one of the following two situation occurs:

a) $v(Y_1)$ divides $v(Z_2)$ and the residue of $Z_2/Y_1^r$ is not in
$k$ with $v(Z_2)=r\cdot v(Y_1)$.

b) $v(Y_1)$ does not divide $v(Z_2)$.

In case a), we make the transformation
$$Z_2=Y_2-\sum_{i=1}^{s_2}a_{2,i}Y_1^i,$$
and apply lemma \ref{lema1} to obtain elements with the same values. We
note these elements by $Y_1,\ldots ,Y_n$ again in order not to
complicate the notation. So, after this procedure, we have a
transcendental element $u_2=\si (Y_2/Y_1+\M)$ over $k$.

In case b) we make the same transformation and go back to the
beginning of the proof.

Anyway this procedure stops, because the value of the variables are
greater or equal than 1.

Then we can suppose that, after a finite number of
transformations, we have $n$ elements $Y_1,\ldots ,Y_n$ such that
$v(Y_i)=v(Y_1)=\al$ and the residue $Y_2/Y_1+\M$ is not in $k$.\qed
\end{pf}

\begin{exmp}
Let $v=\nu_t\circ\Psi$ the discrete valuation of $\CC ((X_1,X_2))\vert\CC$
defined by the embedding
$$\begin{array}{rcl}
\Psi :\CC\lcor X_1,X_2\rcor & \lto & \CC (u)\lcor t\rcor \\
X_1 & \lmapsto & t\\
X_2 & \lmapsto & t+t^3+\sum_{i\geq 1}u^it^{i+3}
\end{array}
$$
with $u$ and $t$ independent variables over $\CC$.

The residue $X_2/X_1+\m =1+\m$, because $v(X_2-X_1)=3>1$. So we have
$$v\left(\frac{X_2-X_1}{X_1^3}\right) =0.$$
The residue 
$$\frac{X_2-X_1}{X_1^3}+\m =1+\m$$ 
too, because $v(X_2-X_1-X_1^3)=4>3$. So we have
$$v\left(\frac{X_2-X_1-X_1^3}{X_1^4}\right) =0.$$
As $\Psi ((X_2-X_1-X_1^3)/X_1^4)=u$ and $u$ is trancendental over $\CC$,
then
$$\frac{X_2-X_1-X_1^3}{X_1^4}+\m\notin \CC$$
and this is a first transcendental residue of $\d$ over $\CC$.

In this situation we can do the transformation
$$\begin{array}{rcl}
\CC\lcor X_1,X_2\rcor & \lto & \CC\lcor Y_1,Y_2\rcor\\
X_1 & \lmapsto & Y_1\\
X_2 & \lmapsto & Y_2Y_1^3+Y_1+Y_1^3
\end{array}
$$
to obtain elements $\{ Y_1,Y_2\}$ such that $\Psi (Y_1)=t$ and $\Psi
(Y_2)=\sum_{i\geq 1}u^it^i$. So $v(Y_2)=v(Y_1)=1$ and the residue
$$\frac{Y_2}{Y_1}+\m =\frac{X_2-X_1-X_1^3}{X_1^4}+\m$$
is not in $\CC$.

In this example the extension of the valuation $v$ to the field $\CC
((Y_1,Y_2))$ is the usual order function. Theorem \ref{thm3} says that, for
$n=2$, we always have this.
\end{exmp}

We end up the section with some specific arguments for the case $n=2$.

The proof of the following lemma is straightforward from (\cite{Bri2},
theorem 2.4):

\begin{lem}
Let $v$ be a discrete valuation of $K_n\vert k$. If $v$ is such that
$v(f_r)=r\al$ for all forms $f_r$ of degree $r$ with respect to the usual
degree, then the group of $v$ is $\al\cdot\ZZ$.
\end{lem}

So we have

\begin{thm}{\label{thm3}}
In the case $n=2$, the extension of the valuation $v$ to the field
$k((Y_1,Y_2))$ is the usual order function.
\end{thm}

\begin{pf}
After a finite number of transformations we are in the situation
of the end of the previous proof. Obviously, if $n=2$,
$k((Y_1,Y_2))\subset\R$ so $v$ can be extended to a valuation $v'$
over $k((Y_1,Y_2))$ such that $\Delta_{v'}=\Delta_v=\Delta_{\v}$.
We denote the extension by $v$ for simplifying. Let $\si :\D\to\R$ a
$k-$section of $\R\to\D$, $u_2=\si (Y_2/Y_1+\M )$, $h\ne 0$ a form
of degree $r$ and $\gamma =Y_2-u_2Y_1$. From the construction
procedure of
$u_2$ we know that $\v (\gamma )> \al$ (remember $\al
=v(Y_1)$). Then
$$h(Y_1,Y_2)=h(Y_1,u_2Y_1+\gamma
)=Y_1^rh(1,u_2)+\gamma ',$$
where $\gamma '$ is such that
$v(\gamma ')>r\al$. As $u_2\notin k$, $u_2$ is
transcendental over $k$, so $h(1,u_2)\ne 0$ and $v(h)=r\al$. By
the previous lemma, the group of $v$ is $\al\cdot\ZZ$, so
$\al =1$ and $v$ is the usual order function. \qed
\end{pf}

\section{The general case}

Let us move to the general case. Assume that $n>2$ and suppose we have
applied the procedure of the lemma \ref{lema6} to find $Y_1,\ldots
,Y_n\in\K$ such that

a) The value of these elements are $\al\in\ZZ$.

b) The residue of $Y_2/Y_1$ is transcendental over $k$.

This section and the next one describe the transformations that we have
to do in order to construct the residue field of $v$.

\begin{rem}
Let $\Delta_2=k(Y_2/Y_1+\m )$ a purely
transcendental extension of $k$ of transcendence degree 1. Let
$\si_2:\De_2\to k(Y_2/Y_1)$ defined by
$$\si_2\(\frac{Y_2}{Y_1}+\m\) =\frac{Y_2}{Y_1}=u_2.$$ We know that
there exists a $k-$section $\si$ which extends $\si_2$ in the sense of
the remark \ref{remark4}.
\end{rem}

\begin{rem}
Let  us suppose that the residue of $Y_3/Y_1$
is algebraic over $\Delta_2$, and let $u_{3,1}$ be its image by
$\si$. Then $v(Y_3-u_{3,1}Y_1)=\al_1>\al$. If $\al$ divides to
$\al_1$ then there exists $u_{3,r}\in\im (\si )$ and $r>1$ such
that $v(Y_3-u_{3,1}Y_1-u_{3,r}Y_1^r)=\al_2>\al_1$. Let us suppose
that $u_{3,r}$ is algebraic over $\Delta_2$ too and $\al$ divides
to $\al_2$. Then we can find ourselves in one of the three
situations shown in the following items.

{\bf (Situation 1)} After a finite
number of transformations, we obtain a value $\al_s$ such that
it is not divided by $\al$. Then we make the transformation
$$Z_3=Y_3-\sum_{j=1}^su_{3,j}Y_1^j,$$ with $u_{3,j}$ algebraic over
$\Delta_2$ for all $j=1,\ldots ,s$. So we have to apply
transformations to find elements with the same values and begin with
all the procedure described in this section.
When this occurs, the values of the elements decrease, so we can suppose that
after a finite number of transformations we have reached a strictly
minimal value. In fact this value should be 1, because we are assuming that
the values group is $\ZZ$. We shall denote these elements
by $Y_1,\ldots ,Y_n$ in order not to complicate the notation. So
we can suppose that this situation will never occur again for any
variable.

\noindent {\bf (Situation 2)} After a finite
number of steps, we have a transcendental residue of
$\Delta_2$. Let us denote this residue by $u_3$. This means
$$Z_3=Y_3-\sum_{j=1}^{s_3}u_{3,j}Y_1^j,$$
where the elements $\{ u_{3,j}\}_{j=1}^{s_3}$ are algebraic over $\Delta_2$
and $u_3=\si (Z_3/Y_1^{\v (Z_3)}+\m )$ is transcendental over
$\Delta_2$. We shall note $\Delta_3 =k(u_2,\{ u_{3,j}\}_{j=1}^{s_3},u_3)$.

In this situation, if $n=3$ we can apply monoidal
transformations to obtain elements with the same values. We will denote
these elements again by $\{ Y_1,Y_2 ,Y_3\}$. The extension of the
valuation $v$ to the field
$L((Y_1,Y_2,Y_3))$ with $L=k(\{u_{3,j}\}_{j=1}^{s_3})$, is the usual
order function, for analogy with the case $n=2$ (theorem 9).

\noindent {\bf (Situation 3)} All the residues
obtained are algebraic elements. Then we take $\Delta_3=\Delta_2(\{
u_{3,j}\}_{j\geq 1})$, an algebraic extension of $\Delta_2$.
\end{rem}


\begin{rem}
Let us suppose that we have repeated the previous construction with each
element $Y_4,\ldots ,Y_{i-1}$, so we have a field
$$\De_{i-1}=k(u_2,\zeta_3,\ldots ,\zeta_{i-1})\subset \si (\D ),$$
where each $\zeta_k$ is: 

- either $\{\{ u_{k,j}\}_{j=1}^{s_k},u_k\}$ if $\{
u_{k,j}\}_{j=1}^{s_k}$ are algebraic over $\De_{k-1}$ and $u_k=\si
((Z_k/Y_1^{\v (Z_k)})+\M )$ is a transcendental element over $\De_{k-1}$
(i.e. situation 2), 

- or $\De_{k-1}\subset\De_{k-1}(\{ u_{k,j}\}_{j\geq 1})$
is an algebraic extension (i.e. situation 3). 

So we have two possible situations concerning variable $Y_i$: 

{\bf 1)} There exists a transformation
$$Z_i=Y_i-\sum_{j=1}^{s_i}u_{i,j}Y_1^j,$$ 
where the elements $u_{i,j}$ are algebraic over $\De_{i-1}$ and $u_i=\si
((Z_i/Y_1^{\v (Z_i)})+\M )$ is a transcendental element over $\De_{i-1}$. So
we have the transcendental extension 
$$\De_{i-1}\subset\De_{i-1}(\{u_{i,j}\}_{j=1}^{s_i},u_i)=\De_i.$$

{\bf 2)} All the elements $u_{i,j}$ we have
constructed are algebraic over $\De_{i-1}$, so we have the
algebraic extension
$$\De_{i-1}\subset\De_{i-1}(\{u_{i,j}\}_{j\geq 1})=\De_i.$$
\end{rem}

\begin{rem}
We have given a procedure to  construct elements $\{ Y_1,\ldots
,Y_n\}$ such that they satisfy these important properties:
\begin{enumerate}
\item After reordering if necessary, we can suppose that the first $m$
elements give us all the
transcendental residues over $k$, i.e. the residue of each
$Y_i/Y_1$ is transcendental over $\De_{i-1}$ with $i=2,\ldots ,m$.
So the rest of variables $Y_{m+1},\ldots ,Y_n$ are such that we
enter in the situation of previous item {\bf 2)}.
\item With the usual notations, the extension
$$\De_m\subset\De_m\(\{u_{i,j}\}_{j\geq 1}\),\ i=m+1,\ldots ,n$$
is algebraic.
\end{enumerate}
\end{rem}


\begin{thm}{\label{theorem9}}
The residue field of $v$ is
$$\De_n=k\( u_2,\{u_{3,j}\}_{j=1}^{s_3},u_3,\ldots
,\{u_{m,j}\}_{j=1}^{s_m},u_m\) \(\{u_{m+1,j}\}_{j\geq 1},\ldots
,\{u_{n,j}\}_{j\geq 1}\) ,$$
and the transcendence degree of $\De_n$ over $k$ is $m-1$.
\end{thm}

\begin{pf}
In this section we have given a construction by writing the
elements $Y_i$ depending on $Y_1$ and some transcendental and
algebraic residues. So we have constructed a map
$$
\begin{array}{rcl}
\vp ':L_n(( \Y )) & \longrightarrow & \Delta_n((t )) \\
Y_1 & \lmapsto & t \\
Y_i & \lmapsto & u_it,\ i=2,\ldots ,m\\
Y_k & \lmapsto & \sum_{j\geq 1}u_{k,j}t^j,\ u_{k,1}\ne 0,\
k=m+1,\ldots ,n.
\end{array}
$$ 
This map is not injective in the general case, but we know that
$v=\nu_t\circ\vp '_{\vert K_n}$. So the residue field of $v$ is
equal to the residue field of $\nu_t$, i.e. $\De_n$. \qed
\end{pf}

A straightforward consequence of this theorem is the following
well-known result

\begin{cor}{\label{cor8}}
The usual order function over $K_n$ has dimension $n-1$, i.e. the
transcendence degree of its residue field over $k$ is $n-1$.
\end{cor}

\begin{pf}
Let $\nu$ be the usual order function over $K_n$. All the residues
$X_i/X_1+{\got m}_{\nu}$ are transcendental over $k(X_2/X_1+{\got
m}_{\nu},\ldots , X_{i-1}/X_1+{\got m}_{\nu})$: if this were not the
case, there would exist
$u_i\in\si (\Delta_{\nu})$ such that $\nu (X_i-u_iX_1)>1$ and $\nu$
would not be an order function. So $\Delta_{\nu}=k(X_2/X_1,\ldots ,X_{n}/X_1)$.
\qed
\end{pf}

\section{Explicit construction of the residue field: an example}

In order to compute explicitly the residue field of a valuation we need
to construct a section $\si :\D\to\R$ as in remark \ref{remark4}. This
procedure is not constructive in general. As in section 1, if the
valuation is given as a composition $v=\nu_t\circ\Psi$, where $\Psi
:k\lcor\X\rcor\to\Delta\lcor t\rcor$ is an injective homomorphism and
$\nu_t$ is the order funcion in $\Delta\lcor t\rcor$, then we can construct
$\si$ using the coefficients $a_{i,j}\in\Delta$ of $\Psi (X_i)=\sum_{j\geq
1}a_{i,j}t^j$.

(Of course, explicit does not mean effective because we are working with the
series $\sum_{j\geq 1}a_{i,j}t^j$ and this input is not finite).

\begin{exmp}{\label{ejemplo}}
Let us consider the embedding
$$
\begin{array}{rcl}
\Psi :  \CC\lcor X_1,X_2,X_3,X_4,X_5\rcor & \lto &
        \Delta\lcor t\rcor \\
X_1 & \lmapsto & t \\
X_2 & \lmapsto & T_2t \\
X_3 & \lmapsto & T_2^2t+T_2t^2+T_3t^3 \\
X_4 & \lmapsto & T_2^3t+T_2^2t^2+T_3t^3+T_4t^4\\
X_5 & \lmapsto & T_2t\sum_{j\geq 1} (T_4^{1/p}t)^j,
\end{array}
$$
with $t$, $T_2$, $T_3$ and $T_4$ variables over $\CC$, $p\in\ZZ$ prime and
$\Delta$ is a field such that $\overline{\CC (T_4)}(T_2,T_3)\subseteq\Delta$.
$\overline{\CC (T_4)}$ is the algebraic closure of $\CC (T_4)$. We are going
to denote  its extension to the quotient fields by $\Psi$. The
composition of this injective homomorphism with the order function in
$t$ gives a discrete valuation of $\CC ((X_1,X_2,X_3,X_4,X_5))\vert\CC$,
$v=\nu_t\circ\Psi$. The residues of $X_i/X_1$ are not in $\CC$ for
$i=2,3,4,5$.

Let us put $u_2=\si (X_2/X_1+\m )$, a transcendental element over
$\CC$. By remark \ref{remark4} we know how to construct $\si$ step by
step, so let take us $u_2=X_2/X_1$ and $\De_2=\CC(u_2)$.

The residue $X_3/X_1+\m$ is algebraic over $\CC (u_2)$, in fact
$$\frac{X_3}{X_1}+\m =\frac{X_2^2}{X_1^2}+\m .$$
So we can take
$u_{3,1}=\si ((X_3/X_1)+\m )=u_2^2$. The value of $X_3-u_{3,1}X_1$
is 2, therefore we have to see if the residue
$$\frac{X_3-u_{3,1}X_1}{X_1^2}+\m$$
is algebraic over $\CC (u_2)$.
We have that $$\frac{X_3-u_{3,1}X_1}{X_1^2}+\m = \frac{X_2}{X_1}+\m
,$$
so it is algebraic and we can take $u_{3,2}=u_2$. Now
$v(X_3-u_{3,1}X_1-u_{3,2}X_1^2)=3$ and we have to check if
$$\frac{X_3-u_{3,1}X_1-u_{3,2}X_1^2}{X_1^3}+\m$$ is algebraic over
$\De_2$. In this case, as
$$\Psi\(\frac{X_3-u_{3,1}X_1-u_{3,2}X_1^2}{X_1^3}+\m\)= T_3 ,$$ this
residue is transcendental. So we take
$$u_3=\si\(\frac{X_3-u_{3,1}X_1-u_{3,2}X_1^2}{X_1^3}+\m\)
=\frac{X_1X_3-X_2^2-X_1^2X_2}{X_1^4} .$$ Let us take $\De_3=\CC
(u_2,u_3)$.

We have to apply this procedure to $X_4$. The
residue $X_4/X_1+\m$ is algebraic over $\De_3$ because
$$\frac{X_4}{X_1}+\m =\frac{X_2^3}{X_1^3}+\m ,$$ 
so we can take $u_{4,1}=\si ((X_4/X_1)+\m )=u_2^3\in\De_3.$
Now $v(X_4-u_{4,1}X_1)=2$, and we have to check what happens with the
residue
$$\frac{X_4-u_{4,1}X_1}{X_1^2}+\m .$$ 
As
$$\frac{X_4-u_{4,1}X_1}{X_1^2}+\m = \frac{X_1^2}{X_2^2}+\m ,$$ 
it holds 
$$u_{4,2}=\si\(\frac{X_4-u_{4,1}X_1}{X_1^2}+\m\) =u_2^2.$$
Clearly $v(X_4-u_{4,1}X_1-u_{4,2}X_1^2)=3$ and
$$\frac{X_4-u_{4,1}X_1-u_{4,2}X_1^2}{X_1^3}+\m
=\frac{X_1X_3-X_2^2-X_1^2X_2}{X_1^4}+\m ,$$ 
therefore
$$u_{4,3}=\si\(\frac{X_4-u_{4,1}X_1-u_{4,2}X_1^2}{X_1^3}+\m\) =u_3.$$
The following residue is transcendental because
$v(X_4-u_{4,1}X_1-u_{4,2}X_1^2-u_{4,3}X_1^3)=4$ and
$$\Psi\(\frac{X_4-u_{4,1}X_1-u_{4,2}X_1^2-u_{4,3}X_1^3}{X_1^4}\)=T_4.$$
Then we can take
$$u_4=\si\(\frac{X_4-u_{4,1}X_1-u_{4,2}X_1^2-u_{4,3}X_1^3}{X_1^4}+\m\)
=$$
$$=\frac{X_1^2X_4-X_3^2-X_1^2X_2^2-X_1^2X_3-X_1X_2^2-X_1^2X_2}{X_1^6}.$$
So $\De_4=\CC (u_2,u_3,u_4)$.

With the variable $X_5$ we obtain the next algebraic residues
$$u_{5,j}=\si\(\frac{X_5-u_{5,1}X_1-\cdots -u_{5,j-1}X_1^{j-1}}{X_1^j}+\m\)
=u_4^{1/p^j}$$
for all $j\geq 1$. So we have $\Delta_5=\CC (u_2,u_3,u_4)(\{
u_4^{1/p^j}\}_{j\geq 1} )$, an algebraic extension of $\Delta_4$.

Then the residue field of $v$ is
$$\d =\CC\(\frac{X_2}{X_1}+\m ,\frac{X_1X_3-X_2^2-X_1^2X_2}{X_1^4}+\m
, \right.$$
$$\left.\frac{X_1^2X_4-X_3^2-X_1^2X_2^2-X_1^2X_3-X_1X_2^2-X_1^2X_2}{X_1^6}+\m\)
 \(\left\{\(\frac{X_2}{X_1}+\m\)^{1/p^j}\right\}_{j\geq 1}\).$$

In this case, by the transformation
$$
\begin{array}{rcl}
X_1 & \lto & Y_1 \\
X_2 & \lto & Y_2 \\
X_3 & \lto & Y_1^2Y_3+u_{3,1}Y_1+u_{3,2}Y_1^2 \\
X_4 & \lto & Y_1^3Y_4+u_{4,1}Y_1+u_{4,2}Y_1^2+u_{4,3}Y_1^3 \\
X_5 & \lto & Y_5,
\end{array}
$$
we can extend the valuation $v$ to a discrete valuation $v'=\nu_t\Psi'$ of
$\CC ((Y_1,Y_2,Y_3,Y_4,Y_5))$, with the injective homomorphism 
$$
\begin{array}{rcl}
\Psi' :  \CC\lcor Y_1,Y_2,Y_3,Y_4,Y_5\rcor & \lto &
        \Delta\lcor t\rcor \\
Y_1 & \lmapsto & t \\
Y_i & \lmapsto & T_it,\ i=2,3,4\\
Y_5 & \lmapsto & \sum_{j\geq 1} (T_4^{1/p}t)^j.
\end{array}
$$
The restriction $v'_{\vert\CC ((Y_1,Y_2,Y_3,Y_4))}$ is the usual
order function. This is not the general case because $\Psi'$ may not be
injective.
\end{exmp}

\section{Rank one discrete valuations and order functions}

We can summarize the constructions of previous sections in the following
theorem wich generalize the results of \cite{Bri2,Br-He}

\begin{thm}{\label{teorema212}}
Let $v$ be a discrete valuation of $K_n\vert k$, then
\begin{enumerate}
\item If the dimension of $v$ is $n-1$, we can embed $k\lcor\X\rcor$
into a ring $L\lcor\Y\rcor$, where $L\subset\si (\D )$ and the
extended valuation of $v$ over the field $L((\Y ))$ is the usual order
function.
\item If the dimension of $v$ is $m-1<n-1$, we can embed $k\lcor\X\rcor$
into a ring $L\lcor\Y\rcor$, where $L\subset\si (\D )$ and the
restriction into $L((\Ym ))$ of the ``extended valuation'' of $v$ over
$L((\Y ))$ is the usual order function.
\end{enumerate}
\end{thm}

\begin{pf} We have the following map:
$$
\begin{array}{rcl}
\vp ':L_n(( \Y )) & \longrightarrow & \Delta_n((t )) \\
Y_1 & \lmapsto & t \\
Y_i & \lmapsto & u_it,\ i=2,\ldots ,m \\
Y_k & \lmapsto & \sum_{j\geq 1}u_{k,j}t^j,\ u_{k,1}\ne 0,\
k=m+1,\ldots ,n,
\end{array}
$$
where $m-1$ is the dimension of $v$. Let us prove the theorem:

\begin{enumerate}

\item In the case $m=n$, $\vp '(Y_i)=u_it$ for all $i=2,\ldots
,n$. Let $\nu_t$ be the usual order funtion over $\Delta_n ((t))$. The
homomorphism $\vp '$ is injective and the valuation $v'=\nu_t\circ\vp '$
of $L((\Y ))$ is the usual order fuction over this field.
Obviously $v'$ extends $v$.

\item If $m<n$ we can consider the elements $W_k=Y_k-\sum_{j\geq
1}u_{k,j}Y_1^j$. Hence we have $L ((\Y))
= L((Y_1,\ldots ,Y_m,W_{m+1},\ldots ,W_n))$. We define the
discrete valuation of rank $n-m+1$ over $L(( \Y ))$:
$$v' (Y_1)=\ldots
=v'(Y_m)= (0,\ldots ,0,1),$$
$$v'(W_{m+1})=(0,\ldots ,1,0),\ldots
,v'(W_n)=(1,0,\ldots ,0).$$
The restriction of this valuation to $K_n$
is a rank one discrete valuation, because the value of any element is
in $0\times\cdots\times 0\times\ZZ$. In fact $v'(f)=(0,\ldots
,0,v(f))$ for all $f\in K_n$, so $v'$ ``extends'' $v$ in this
sense. Obviously $v'_{\vert L((\Ym ))}$ is the usual order
function. We want note that this ideal $(W_{m+1},\ldots ,W_{n})$ is
the {\em implicit ideal of $v$} that appears in some works of
M. Spivakovsky (\cite{Sp}).
\end{enumerate}
\qed
\end{pf}

For the case of valuations of dimension $n-1$, we can combine
corollary \ref{cor8} and assertion 1 of the previous theorem:

\begin{cor}{\label{theorem8}}
Let $v$ be a discrete valuation of $K_n\vert k$. The following
conditions are equivalent:

1) The transcendence degree of $\D$ over $k$ is $n-1$.

2) There exists a finite sequence of monoidal transformations and
   coordinates changes which take $v$ into an order function.
\end{cor}

\begin{exmp}
Let us consider the homomorphism
$$
\begin{array}{rcl}
\Psi :  \CC\lcor X_1,X_2,X_3,X_4,X_5\rcor & \lto &
        \Delta\lcor t\rcor \\
X_1 & \lmapsto & t \\
X_2 & \lmapsto & T_2t \\
X_3 & \lmapsto & T_2^2t+T_2t^2+T_3t^3 \\
X_4 & \lmapsto & T_2^3t+T_2^2t^2+T_3t^3+T_4t^4\\
X_5 & \lmapsto & T_2t\left(\sum_{j\geq 1} a_j(T_4t)^j\right) ,
\end{array}
$$
with $a_j\in\CC$ such that $\Psi$ is injective (we can take $\sum_{j\geq 1}
a_j(T_4t)^j= e^{T_4t}-1$). Then the residue field of this valuation (see
example \ref{ejemplo}) is
$$\d =\CC\(\frac{X_2}{X_1}+\m ,\frac{X_1X_3-X_2^2-X_1^2X_2}{X_1^4}+\m
, \right.$$
$$\left.\frac{X_1^2X_4-X_3^2-X_1^2X_2^2-X_1^2X_3-X_1X_2^2-X_1^2X_2}{X_1^6}+\m\)
.$$
By the transformation (see example \ref{ejemplo})
$$
\begin{array}{rcl}
X_1 & \lto & Y_1 \\
X_2 & \lto & Y_2 \\
X_3 & \lto & Y_1^2Y_3+u_{3,1}Y_1+u_{3,2}Y_1^2 \\
X_4 & \lto & Y_1^3Y_4+u_{4,1}Y_1+u_{4,2}Y_1^2+u_{4,3}Y_1^3 \\
X_5 & \lto & Y_2Y_5,
\end{array}
$$
we obtain a new field $\CC ((Y_1,Y_2,Y_3,Y_4,Y_5))$, but we can not extend
$v$ to this field because the homomorphism
$$
\begin{array}{rcl}
\Psi' :  \CC\lcor Y_1,Y_2,Y_3,Y_4,Y_5\rcor & \lto &
        \Delta\lcor t\rcor \\
Y_1 & \lmapsto & t \\
Y_i & \lmapsto & T_it,\ i=2,3,4 \\
Y_5 & \lmapsto & \sum_{j\geq 1} a_j(T_4t)^j
\end{array}
$$
is not injective. Then let us take $W_5=Y_5-\sum_{j\geq 1} a_j(Y_4)^j$
(because we can consider $T_4Y_1=Y_4$). Then $\CC
((Y_1,Y_2,Y_3,Y_4,Y_5))=\CC ((Y_1,Y_2,Y_3,Y_4,W_5))$ and the discrete
valuation of rank 2 defined by $v'(Y_i)=(0,1)$ for $i=1,\ldots ,4$ and
$v'(W_5)=(1,0)$ is such that for all $f\in\CC ((X_1,X_2,X_3,X_4,X_5))$ we
have $v'(f)=(0,v(f))$ and $v'_{\vert\CC ((Y_1,Y_2,Y_3,Y_4))}$ is the usual
order function.
\end{exmp}


\begin{thebibliography}{1}


\bibitem{Bri2}
E.~Briales, {Constructive theory of valuations}, \emph{Comm. Algebra}
\textbf{17} (1989), no.~5, 1161--1177.

\bibitem{Br-He}
E.~Briales and F.J. Herrera, {Construcci\'{o}n expl\'{\i}cita de las
  valoraciones de un anillo de series formales en dos variables}, \emph{Actas X
  Jornadas Hispano-Lusas} (Murcia, 1985) vol.~II, pp.~1--10.

\bibitem{Krull}
W.~Krull, {Allgemeine Bewerstungstheorie}, \emph{J. Reine
  Angew. Math.} \textbf{167} (1931), 160--196.


\bibitem{Ser1}
J.~P. Serre, \emph{Corps locaux}, Hermann, Paris, 1968.

\bibitem{Sp}
M.~Spivakovsky, \emph{Resolution of singularities}, preprint, Dept. Of
Math., University of Toronto, August 1994.

\bibitem{ZSII}
O.~Zariski and P.~Samuel, \emph{Commutative algebra}, Graduate Texts in
  Mathematics, vol.~II, Springer-Verlag, New York-Heidelberg-Berlin, 1975.

\end{thebibliography}

\providecommand{\bysame}{\leavevmode\hbox to3em{\hrulefill}\thinspace}

\end{document}